\documentclass{amsart}
\usepackage{amssymb}
\newtheorem{lemma}{Lemme}[section]
\newtheorem{theorem}[lemma]{Theorem}
\newtheorem{proposition}[lemma]{Proposition}
\newtheorem{corollary}[lemma]{Corollary}

\theoremstyle{definition}

\theoremstyle{remark}

\newtheorem{remark}[lemma]{Remark}

\def\SZ{\mathbb Z}    %%% Z
\def\SQ{\mathbb Q}    %%% Q
\def\pp{\mathfrak p}    %%% p
\def\qq{\mathfrak q}    %%% q
\def\int{\mbox{\rm{Int}}}             %%% Int

\def\Oo{\mathcal{O}}

\def\Ii{\mathcal{I}}

\def\Gal{\textrm{Gal}}

\def\be{\begin{equation}}
\def\ee{\end{equation}}

\begin{document}

\title[On Abhyankar's lemma]{On Abhyankar's lemma\\about ramification indices}

\author{Jean-Luc Chabert}
\address{LAMFA CNRS-UMR 7352\\Universit\'e de Picardie\\80039 Amiens, France}
\email{jean-luc.chabert@u-picardie.fr}
\author{Emmanuel Halberstadt}
%\address{Universit\'e Pierre et Marie Curie-Paris VI}

\begin{abstract}
We provide a simple proof of the fact that the ramification index of the {\em compositum} of two finite extensions of local fields is equal to the least common multiple of the ramification indices when at least one of the extensions is tamely ramified.
\end{abstract}

%%%%%%%%%%%%%%%%%%%%%%%%%%%%%%%%%%%%%%%%
\maketitle

%%%%%%%%%%%%%%%%%%%%%%

\section{Introduction}

Let $L/K$ be an extension of number fields and assume that $L$ is the {\em compositum} of two sub-extensions $K_1$ and $K_2$. Let $\qq$ be a prime ideal of $L$ and let $\pp=\qq\cap K$, $\pp_1=\qq\cap K_1$, $\pp_2=\qq\cap K_2$. We denote by $e(\qq/\pp)$ the ramification index of $\qq$ in the extension $L/K$. Then, by multiplicativity of the ramification indices, that is,
$$e(\qq/\pp)=e(\qq/\pp_i)\times e(\pp_i/\pp) \quad (i=1,2)\,,$$
we obviously have:
\be\label{eq:1} \mathrm{lcm}\,\{e(\pp_1/\pp),e(\pp_2/\pp)\}\;\big\vert\; e(\qq/\pp) \,.\ee
On the other hand, if one of the extensions $K_i/K$ is normal, for instance if $K_1/K$ is normal, the extension $L/K_2$ is normal and the following morphism is injective:
$$\rho:\sigma\in\Gal(L/K_2)\mapsto\sigma_{\vert K_1}\in\Gal(K_1/K) \,.$$
Recall that, as the residue fields are perfect, the ramification index $e(\qq/\pp_2)$ is equal to the order of the inertia group $\Ii_{\qq}(L/K_2)$ of $\qq$ in the extension $L/K_2$, that is,
$$\Ii_{\qq}(L/K_2)=\{\sigma\in\Gal(L/K_2)\mid \forall x\in\Oo_L\;\;\sigma(x)-x\in\qq\}\,.$$
Now, the image by $\rho$ of $\Ii_{\qq}(L/K_2)$ is clearly contained in the inertia group $\Ii_{\pp_1}(K_1/K).$ Thus, $e(\qq/\pp_2)$ divides $e(\pp_1/\pp)$, and hence,
\be\label{eq:2} e(\qq/\pp) \; \big\vert \; e(\pp_1/\pp)\times e(\pp_2/\pp)\,.\ee 
Formula (\ref{eq:2}) may be false in general (see Remark~\ref{rem} below).

\medskip
There is another well known result about ramification indices of {\em composita}, namely 
Abhyankar's lemma. This result is generally known in the following form:

\begin{proposition} $($\emph{Narkiewicz} \cite[p. 229]{bib:nark}$)$\label{th:1}
If $K,K_1,K_2$ are local fields such that $K_1/K$ is tame, $K_2/K$ is finite and $e(K_1/K)$ divides $e(K_2/K)$, then $K_1K_2/K_2$ is unramified.
\end{proposition}

Roughly speaking, one may kill tame ramification by taking an extension of the base field (see also \cite[p. 279]{bib:SGA}).
In fact, one finds in \cite{bib:stich} a stronger formulation, but it is stated only for function fields:

\begin{proposition} $($\emph{Stichtenoth} \cite[Th. 3.9.1]{bib:stich}$)$ \label{th:2}
Let $L/K$	be a finite separable extension of function fields. Suppose that $L=K_1K_2$ is the compositum of two intermediate fields $K\subseteq K_1,K_2\subseteq L$. Let $Q$ be a place of $L$ extension of a place $P$ of $K$ and set $P_i:=Q\cap K_i$ for $i=1,2$. Assume that at least one of the extensions $P_1 / P$ or $P_2 / P$ is tame. Then
$$e(Q / P)=\mathrm{lcm}\{e(P_1 / P),e(P_2 / P)\}\,.$$
\end{proposition}

Since we did not find in the literature such a statement with respect to number fields (although it probably exists somewhere hidden under an indirect formulation), we provide here a simple proof of this generalized result (proof which in some sense is close to that of Proposition~\ref{th:2}).

%%%%%%%%%%%%

\section{Abhyankar's lemma}

\begin{theorem}\label{th}
Let $A$ be a Dedekind domain with quotient field $K$. Let $L/K$ be a finite separable extension of fields. Assume that $L$ is the {\em compositum} of two subfields $K_1$ and $K_2$. Denote by $B$, $B_1$ and $B_2$ the integral closures of $A$ in $L$, $K_1$ and $K_2$. Let $\pp$ be a maximal ideal of $A$ whose residue field $A/\pp$ is perfect with characteristic $p\;($which may be $0)$. Finally, let $\qq$ be a maximal ideal of $B$ lying over $\pp$ and let $\pp_i=\qq\cap A_i$ for $i=1,2$. 

If at least one of the extensions $K_i/K$ is tamely ramified in $\pp_i\;($that is, if one of the integers $e(\pp_i/\pp)$ is not divisible by $p)$, then one has the equality:
\be e(\qq/\pp)=\mathrm{lcm}\,\{\,e(\pp_1/\pp),e(\pp_2/\pp)\,\}\,.  \ee
\end{theorem}

Note that, if the characteristic of $A/\pp$ is 0, then the ramification is always tame. Of course, Propositions~\ref{th:1} and \ref{th:2} are consequences of Theorem~\ref{th}.

\begin{proof}
Let $L'$ be the normal closure of $L$ over $K$, let $B'$ be the integral closure of $B$ in $L'$ and let $\qq'$ be a maximal ideal of $B'$ lying over $\qq$. One knows that $e(\qq'/\pp)=\vert G_0\vert$ where $G_0$ denotes the inertia group of $\qq'$ in the extension $L'/K$. Moreover, denoting by $\pi\in B'$ a uniformizer with respect to $\qq'$, we have a group homomorphism:
$$s\in G_0 \mapsto s(\pi)/\pi \textrm{ mod } \qq' \in (B'/\qq')^*$$
whose kernel is the subgroup:
$$G_1=\{s\in\Gal(L'/K)\mid \forall x\in B'\; s(x)-x\in\qq'^2\}\,.$$ 
Thus, $G_1 \lhd G_0$. The injectivity of the morphism $G_0/G_1\to (B'/\qq')^*$
shows that the group $G_0/G_1$ is cyclic	 and that its order is prime to the characteristic $p$ of $B'/\qq'$. Moreover, one knows also that, if $p=0$, then $G_1=\{1\}$  and, if $p>0$, then $G_1$ is a $p$-group (cf., for instance, \cite[IV, \S 2]{bib:serre}). Finally, $G_0$ is a semidirect product of a cyclic group of order prime to $p$ with the normal $p$-group $G_1.$ 

If $E$ is a field between $K$ and $L'$, the analogs of the groups $G_j$ for $j=0,1$ with respect to $\qq'$ in the extension $L'/E$ are clearly the groups $G_j\cap \textrm{Gal}(L'/E)$. Let $\Gamma_1=\textrm{Gal}(L'/K_1)$, $\Gamma_2=\textrm{Gal}(L'/K_2)$ and $\Gamma=\textrm{Gal}(L'/L)$. One has $\Gamma=\Gamma_1\cap \Gamma_2$ since $L=K_1K_2$. Then, by multiplicativity, one has:
$$e(\qq/\pp)=\frac{e(\qq'/\pp)}{e(\qq'/\qq)}= \frac{\vert G_0\vert}{\vert G_0\cap\Gamma\vert} \;\textrm{ and } \; e(\pp_i/\pp)=\frac{e(\qq'/\pp)}{e(\qq'/\pp_i)}= \frac{\vert G_0\vert}{\vert G_0\cap\Gamma_i\vert}\; (i=1,2)\,.$$

Let $m=\mathrm{lcm}\,\{\,e(\pp_1/\pp),e(\pp_2/\pp)\,\}$, then $m=\mathrm{lcm}\,\{\,\frac{\vert G_0\vert}{\vert G_0\cap\Gamma_1\vert},\frac{\vert G_0\vert}{\vert G_0\cap\Gamma_2\vert}\,\}$. Now, let $d=\mathrm{gcd}\,\{\,\vert G_0\cap\Gamma_1\vert,\vert G_0\cap\Gamma_2\vert\,\}$, then
 $m\times d = \vert G_0\vert $. Finally,
$\frac{e(\qq/\pp)}{m}=\frac{d}{\vert G_0\cap\Gamma\vert}\,.$
Thus, we have to prove that $d=\vert G_0\cap\Gamma\vert$. 
\medskip

\noindent{\em $1^{st}$ case}. If $p=0$, $G_0$ is cyclic. Since in a cyclic group the order of a subgroup which is the intersection of two subgroups is the gcd of the orders of these two subgroups, it follows from the equality $\Gamma=\Gamma_1\cap\Gamma_2$ that $\vert G_0\cap\Gamma\vert =d$. 

\medskip

\noindent{\em $2^{nd}$ case}. If $p\not=0$, by hypothesis $p$ does not divide one of the integers $e(\pp_i/\pp)=\frac{\vert G_0\vert}{\vert G_0\cap\Gamma_i\vert}$. Assume that $p$ does not divide $\frac{\vert G_0\vert}{\vert G_0\cap\Gamma_1\vert}$. As $G_1$ is the only $p$-Sylow subgroup of $G_0$, $G_1$ is contained in $\Gamma_1$, and hence, $G_1 = G_1\cap \Gamma_1$ and $G_1\cap\Gamma=G_1\cap\Gamma_2$. 
Thus, trivially, we have: 
\be\label{eq:4a} \vert G_1\cap \Gamma\vert =\textrm{gcd}\{\,\vert G_1\cap\Gamma_1\vert, \vert G_1\cap\Gamma_2\vert\,\}\,.\ee
Moreover, let $\pi:G_0\to G_0/G_1$ be the canonical morphism. Clearly, we have:
$$\vert G_0\cap\Gamma\vert=\vert\pi(G_0\cap\Gamma)\vert \times\vert G_1\cap \Gamma\vert \textrm{ and }\vert G_0\cap\Gamma_i\vert=\vert\pi(G_0\cap\Gamma_i)\vert\times \vert G_1\cap \Gamma_i\vert \;(i=1,2)\,.$$ The containment  $\pi(G_0\cap\Gamma)\subseteq \pi(G_0\cap\Gamma_1)\cap\pi(G_0\cap\Gamma_2)$ is obvious. Let us prove the reverse inclusion. Let $x_i\in G_0\cap\Gamma_i$ $(i=1,2)$ such that $\pi(x_1)=\pi(x_2)$, then $x_2=g_1x_1$ for some $g_1\in G_1\subseteq \Gamma_1$, and hence, $x_2\in G_0\cap\Gamma$. From the equality $\pi(G_0\cap\Gamma) = \pi(G_0\cap\Gamma_1) \cap\pi(G_0\cap\Gamma_2)$ and the fact that the group $G_0/G_1$ is cyclic, we deduce:
\be\label{eq:5a} \vert\pi(G_0\cap\Gamma)\vert =
 \textrm{gcd}\{\, \vert \pi(G_0\cap\Gamma_1)\vert,\vert \pi(G_0\cap\Gamma_2)\vert \,\}\,.\ee 
 Since $\vert G_0\vert /\vert G_1\vert$ and $\vert G_1\vert$ are coprime it follows by multiplicativity from (\ref{eq:4a}) and (\ref{eq:5a}) that:
$$\vert G_0\cap\Gamma\vert =
 \textrm{gcd}\{ \vert G_0\cap\Gamma_1\vert ,\vert G_0\cap\Gamma_2\vert  \}=d\,.$$
\end{proof}

\begin{corollary}
With the same notations as in Theorem~\ref{th}, if $e(\pp_1/\pp)$ and $e(\pp_2,\pp)$ are coprime, then $$e(\qq/\pp)=e(\pp_1/\pp)\times e(\pp_2/\pp)\,.  $$
\end{corollary}

\begin{remark}\label{rem}
When $K_1/K$ and $K_2/K$ are both widely ramified, not only $e(\qq/\pp)$ may be stricly greater than the least common multiple of $e(\pp_1/\pp)$ and $e(\pp_2/\pp)$, but it may happen that it does not divide the product $e(\pp_1/\pp)\times    e(\pp_2/\pp)$. 

For instance, let $K=\SQ$, $K_0=\SQ(j)$, $K_1=\SQ(\sqrt[3]{3})$, $K_2=\SQ(j\sqrt[3]{3})$, and $L=\SQ(j,\sqrt[3]{3})=K_0K_1=K_1K_2$ where $j=\exp(2i\pi/3)$. Then, if $\qq$ is a (in fact, the) prime ideal of $L$ lying over $\pp=3\SZ$, then we have $e(\pp_0/\pp)=2, e(\pp_1/\pp)=e(\pp_2/\pp)=3$. It follows from formula (\ref{eq:1}) that $e(\qq/\pp)=6$ which does not divide $e(\pp_1/\pp)\times e(\pp_2/\pp)=9$. 
\end{remark}

%%%%%%%%%%%%%%%%
%%%%%%%%%%%%%%%%%%%

%%%%%%%%%%%%%%%%%%%%%%%%%%%%%%%%%%%%%%%%%%

%%%%%%%%%%%%%%%%

\end{document}